\documentclass[12pt]{article}
\usepackage{graphicx}
\scrollmode
\usepackage{amsfonts}
\usepackage{amsmath}
\usepackage{amssymb}
\begin{document}
\title{The Full Brownian Web as Scaling Limit of Stochastic Flows}
\author{L.~R.~G.~Fontes\and C.~M.~Newman}
\date{}
\maketitle

\newtheorem{defin}{Definition}[section]
\newtheorem{Prop}{Proposition}
\newtheorem{teo}{Theorem}[section]
\newtheorem{ml}{Main Lemma}
\newtheorem{con}{Conjecture}
\newtheorem{cond}{Condition}
\newtheorem{prop}[teo]{Proposition}
\newtheorem{lem}[teo]{Lemma}
\newtheorem{rmk}[teo]{Remark}
\newtheorem{cor}{Corollary}[section]
\renewcommand{\theequation}{\thesection .\arabic{equation}}
\newcommand{\beq}{\begin{equation}}
\newcommand{\eeq}{\end{equation}}
\newcommand{\beqn}{\begin{eqnarray}}
\newcommand{\beqnn}{\begin{eqnarray*}}
\newcommand{\eeqn}{\end{eqnarray}}
\newcommand{\eeqnn}{\end{eqnarray*}}
\newcommand{\bprop}{\begin{prop}}
\newcommand{\eprop}{\end{prop}}
\newcommand{\bteo}{\begin{teo}}
\newcommand{\bcor}{\begin{cor}}
\newcommand{\ecor}{\end{cor}}
\newcommand{\bcon}{\begin{con}}
\newcommand{\econ}{\end{con}}
\newcommand{\bcond}{\begin{cond}}
\newcommand{\econd}{\end{cond}}
\newcommand{\eteo}{\end{teo}}
\newcommand{\brm}{\begin{rmk}}
\newcommand{\erm}{\end{rmk}}
\newcommand{\blem}{\begin{lem}}
\newcommand{\elem}{\end{lem}}
\newcommand{\ben}{\begin{enumerate}}
\newcommand{\een}{\end{enumerate}}
\newcommand{\bei}{\begin{itemize}}
\newcommand{\eei}{\end{itemize}}
\newcommand{\bdf}{\begin{defin}}
\newcommand{\edf}{\end{defin}}

\newcommand{\nn}{\nonumber}
\renewcommand{\=}{&=&}
\renewcommand{\>}{&>&}
\renewcommand{\le}{\leq}
\newcommand{\+}{&+&}
\newcommand{\fr}{\frac}
\renewcommand{\r}{{\mathbb R}}
\newcommand{\br}{\bar{\mathbb R}}
\newcommand{\Z}{{\mathbb Z}}
\newcommand{\z}{{\mathbb Z}}
\newcommand{\zd}{\z^d}
\newcommand{\zz}{{\mathbb Z}}
\newcommand{\R}{{\mathbb R}}
\newcommand{\tw}{\tilde{\cal W}}
\newcommand{\E}{{\mathbb E}}
\newcommand{\C}{{\mathbb C}}
\renewcommand{\P}{{\mathbb P}}
\newcommand{\N}{{\mathbb N}}
\newcommand{\var}{{\mathbb V}}
\renewcommand{\S}{{\cal S}}
\newcommand{\T}{{\cal T}}
\newcommand{\U}{{\cal U}}
\newcommand{\V}{{\cal V}}
\newcommand{\W}{{\cal W}}
\newcommand{\X}{{\cal X}}
\newcommand{\Y}{{\cal Y}}
\newcommand{\cm}{{\cal M}}
\newcommand{\cp}{{\cal P}}
\newcommand{\h}{{\cal H}}
\newcommand{\f}{{\cal F}}
\newcommand{\cd}{{\cal D}}
\newcommand{\xt}{X_t}
\renewcommand{\ge}{g^{(\epsilon)}}
\newcommand{\xe}{y^{(\epsilon)}}
\newcommand{\ye}{y^{(\epsilon)}}
\newcommand{\bx}{{\bar y}}
\newcommand{\by}{{\bar y}}
\newcommand{\bxe}{{\bar y}^{(\epsilon)}}
\newcommand{\bye}{{\bar y}^{(\epsilon)}}
\newcommand{\bwe}{{\bar w}^{(\epsilon)}}
\newcommand{\bxz}{{\bar y}}
\newcommand{\bwz}{{\bar w}}
\newcommand{\we}{w^{(\epsilon)}}
\newcommand{\Xe}{Y^{(\epsilon)}}
\newcommand{\Ze}{Z^{(\epsilon)}}
\newcommand{\Ye}{Y^{(\epsilon)}}
\newcommand{\ydo}{Y^{(\d)}_{y_0(\d),s_0(\d)}}
\newcommand{\yo}{Y_{y_0,s_0}}
\newcommand{\tye}{{\tilde Y}^{(\epsilon)}}
\newcommand{\hy}{{\hat Y}}
\newcommand{\ve}{V^{(\epsilon)}}
\newcommand{\Ne}{N^{(\epsilon)}}
\newcommand{\ce}{c^{(\epsilon)}}
\newcommand{\cle}{c^{(\l\epsilon)}}
\newcommand{\xet}{Y^{(\epsilon)}_t}
\newcommand{\hxt}{\hat X_t}
\newcommand{\btn}{\bar\tau_n}
\newcommand{\ct}{{\cal T}}
\newcommand{\rn}{{\cal R}_n}
\newcommand{\nt}{{N}_t}
\newcommand{\lnk}{{\cal L}_{n,k}}
\newcommand{\cl}{{\cal L}}
\newcommand{\bw}{\bar{\cal W}}
\newcommand{\tc}{\tilde{\cal C}_b}
\newcommand{\hxtt}{\hat X_{\ct}}
\newcommand{\txnt}{\tilde X_{\nt}}
\newcommand{\xs}{X_s}
\newcommand{\xn}{\tilde X_n}
\newcommand{\tx}{\tilde X}
\newcommand{\hx}{\hat X}
\newcommand{\txi}{\tilde X_i}
\newcommand{\txij}{\tilde X_{i_j}}
\newcommand{\taxi}{\tau_{\txi}}
\newcommand{\txn}{\tilde X_N}
\newcommand{\xk}{X_K}
\newcommand{\ts}{\tilde S}
\newcommand{\tl}{\tilde\l}
\newcommand{\tg}{\tilde g}
\newcommand{\im}{I^-}
\newcommand{\ip}{I^+}
\newcommand{\hal}{H_\a}
\newcommand{\ba}{B_\a}

\renewcommand{\a}{\alpha}
\renewcommand{\b}{\beta}
\newcommand{\g}{\gamma}
\newcommand{\G}{\Gamma}
\renewcommand{\L}{\Lambda}
\renewcommand{\d}{\delta}
\newcommand{\D}{\Delta}
\newcommand{\e}{\epsilon}
\newcommand{\fes}{\phi^{(\epsilon)}_s}
\newcommand{\fet}{\phi^{(\epsilon)}_t}
\newcommand{\fe}{\phi^{(\epsilon)}}
\newcommand{\pset}{\psi^{(\epsilon)}_t}
\newcommand{\pse}{\psi^{(\epsilon)}}
\renewcommand{\l}{\lambda}
\newcommand{\me}{\mu^{(\epsilon)}}
\newcommand{\re}{\rho^{(\epsilon)}}
\newcommand{\tre}{\tilde{\rho}^{(\epsilon)}}
\newcommand{\nue}{\nu^{(\epsilon)}}
\newcommand{\mbe}{{\bar\mu}^{(\epsilon)}}
\newcommand{\rbe}{{\bar\rho}^{(\epsilon)}}
\newcommand{\mb}{{\bar\mu}}
\newcommand{\rb}{{\bar\rho}}
\newcommand{\mbz}{{\bar\mu}}
\newcommand{\s}{\sigma}
\renewcommand{\o}{\Pi}
\newcommand{\om}{\omega}
\newcommand{\tio}{\tilde\o}
\renewcommand{\sl}{\sigma'}
\newcommand{\si}{\s(i)}
\newcommand{\sit}{\s_t(i)}
\newcommand{\ei}{\eta(i)}
\newcommand{\eit}{\eta_t(i)}
\newcommand{\eot}{\eta_t(0)}
\newcommand{\sil}{\s'_i}
\newcommand{\sj}{\s(j)}
\newcommand{\st}{\s_t}
\newcommand{\so}{\s_0}
\newcommand{\xii}{\xi_i}
\newcommand{\xij}{\xi_j}
\newcommand{\xio}{\xi_0}
\newcommand{\ti}{\tau_i}
\newcommand{\te}{\tau^{(\epsilon)}}
\newcommand{\bt}{\bar\tau}
\newcommand{\tti}{\tilde\tau_i}
\newcommand{\tto}{\tilde\tau_0}
\newcommand{\tei}{T_i}
\newcommand{\ttei}{\tilde T_i}
\newcommand{\tes}{T_S}
\newcommand{\tao}{\tau_0}

\renewcommand{\t}{\tilde t}

\newcommand{\da}{\downarrow}
\newcommand{\ua}{\uparrow}
\newcommand{\ar}{\rightarrow}
\newcommand{\lar}{\leftrightarrow}
\newcommand{\va}{\stackrel{v}{\rightarrow}}
\newcommand{\ppa}{\stackrel{pp}{\rightarrow}}
\newcommand{\dw}{\stackrel{w}{\Rightarrow}}
\newcommand{\Va}{\stackrel{v}{\Rightarrow}}
\newcommand{\Ppa}{\stackrel{pp}{\Rightarrow}}
\newcommand{\la}{\langle}
\newcommand{\ra}{\rangle}
\newcommand{\ep}{\vspace{.5cm}}
\newcommand\sqr{\vcenter{
\hrule height.1mm
\hbox{\vrule width.1mm height2.2mm\kern2.18mm\vrule width.1mm}
\hrule height.1mm}}        % This is a slimmer sqr.

\newcommand{\stack}[2]{\genfrac{}{}{0pt}{3}{#1}{#2}}

%%%%%%%%%%%%%%%%%%%%%%%%%%%%%%%%%%%%%%%%%%%%%%%%%%%%%%%%%%%%%%%%%%%%%%%%%%%%%%%%%
%%%%%%%%%%%%%%%%%%%%% ABSTRACT %%%%%%%%%%%%%%%%%%%%%%%%%%%%%%%%%%%%%%%%%%%%%%%%%%%%%%
%%%%%%%%%%%%%%%%%%%%%%%%%%%%%%%%%%%%%%%%%%%%%%%%%%%%%%%%%%%%%%%%%%%%%%%%%%%%%%%%%

\begin{abstract}
In this paper we construct an object which we call the full
Brownian web (FBW) and prove that the collection of all space-time
trajectories of a class of one-dimensional stochastic flows
converges weakly, under diffusive rescaling, to the FBW. The
(forward) paths of the FBW include the coalescing Brownian motions
of the ordinary Brownian web along with bifurcating paths.
Convergence of rescaled stochastic flows to the FBW follows from
general characterization and convergence theorems that we present
here combined with earlier results of Piterbarg.

\end{abstract}


\smallskip

                    
%%%%%%%%%%%%%%%%%%%%%%%%%%%%%%%%%%%%%%%%%%%%%%%%%%%%%%%%%%%%%%%%%%%%%%%%%%%%%%%
%%%%%%%%%%%%%%%% INTRODUCTION %%%%%%%%%%%%%%%%%%%%%%%%%%%%%%%%%%%%%%%%%%%%%%%%%%%%%%
%%%%%%%%%%%%%%%%%%%%%%%%%%%%%%%%%%%%%%%%%%%%%%%%%%%%%%%%%%%%%%%%%%%%%%%%%%%%%%%

\section{Introduction}

\setcounter{equation}{0} \label{sec:intro}

In studying certain one-dimensional stochastic flows~\cite{kn:H},
Piterbarg~\cite{kn:P1,kn:P} show\-ed that on the one hand their
(diffusive) scaling limits led to coalescing Brownian motions
while on the other hand there were necessarily regions of space
time where expansions of trajectories occur rather than the
contractions corresponding to coalescing. Regarding a stochastic
flow as the collection of all its space-time trajectories (or
paths with time $t\in(-\infty,\infty)$), Piterbarg's results
suggested to us that the scaling limit should be expressible as a
collection of paths in space-time related to the Brownian
web~\cite{kn:A2,kn:TW,kn:FINR,kn:FINR1} but, unlike the Brownian
web, including bifurcation as well as coalescence.

One main result of this paper (see Theorem~\ref{teo:sf}) is a
proof that the scaling limit of the collection of all stochastic
flow paths is indeed such an object, which we call the full
Brownian web (FBW). Of course, to prove such a theorem, it is
useful to have first defined the (putative) limiting object and
this is our other main result, the construction (see
Section~\ref{sec:mbw}) of the FBW.

The other crucial results of this paper are a characterization,
given in Section~\ref{sec:char}, of the FBW and a general
convergence theorem (see Theorem~\ref{teo:conv}) based on this
characterization. Theorem~\ref{teo:sf} giving weak convergence of
the rescaled stochastic flow to the FBW  is then an immediate
consequence of our general convergence theorem combined with
Piterbarg's earlier results~\cite{kn:P1,kn:P}. We note that we
only need Piterbarg's result (Theorem 6 in~\cite{kn:P}) concerning
convergence of finite dimensional distributions of forward paths
to coalescing Brownian motions rather than his stronger result
(see Theorem 1 in~\cite{kn:P}) about weak convergence of the
rescaled forward stochastic flow in the Skorohod topology. We
expect that Piterbarg's stronger convergence can be obtained as a
consequence of Theorem~\ref{teo:sf} below, but have not examined
this in detail. We conclude the introduction with more details
about the FBW.

The FBW is a collection of continuous noncrossing paths from $\R$
to $\R$ (which may nonetheless touch, coalesce with or bifurcate
from each other) with the following properties. From every
deterministic space-time point $(x,t)$, there is almost surely a
unique (doubly infinite; i.e., defined for $-\infty<t<\infty$)
path through that point. We will use the term {\em semipath} to
denote both the forward in time semi-infinite path starting at
$(x,t)$ and defined for times $s\geq t$ and the backward in time
semi-infinite path ending at $(x,t)$. The semipath from $(x,t)$ is
distributed as a Brownian path (with unit diffusion coefficient)
starting at that point, and the {\em backward} semipath from
$(x,t)$ is distributed as a backward Brownian path (with unit
diffusion coefficient) from that point. Furthermore, for every
deterministic finite set of space-time points, the forward
semipaths starting from those points are distributed as coalescing
Brownian paths (with unit diffusion coefficient) starting from
those points; an analogous statement holds for the backward
semipaths.

As already mentioned, the FBW arises as the scaling limit of
certain stochastic flows of homeomorphisms~\cite{kn:H,kn:P}, as
defined in Section~\ref{sec:sf}, and can be constructed as a
functional of the (double) Brownian
web~\cite{kn:FINR,kn:FINR1,kn:FINR2} (see also~\cite{kn:TW}). The
structure of the FBW sheds some light on and also raises some
questions about the nature of contracting and expanding
(space-time) regions for these stochastic flows beyond Piterbarg's
analysis in~\cite{kn:P}. Here are three examples.
\begin{enumerate}
    \item According to Proposition~\ref{prop:c2} below (see also
    Remark~\ref{rmk:splice}), every nontrivial path $\gamma$ in the
    FBW consists of a backward Brownian path and a forward Brownian
    path spliced together at some $(x^*,t^*)$. The backward path
    should be regarded as passing through a space-time region of
    expansion, since every point $(x,t)$ with $t<t^*$ along it has
    the property that an FBW path through $(x_1,t)$ for $x_1<x$
    does not coalesce with one through $(x_2,t)$ for $x_2>x$ until
    some time after $t^*$, no matter how close $x_1$ and $x_2$ are
    to $x$ and each other. On the other hand the forward path
    should be regarded as passing through an essentially
    contractive region since for {\em generic} points $(x,t)$ with
    $t>t^*$ along it (generic here means that $(m_{\mbox{{\scriptsize
    \emph{in}}}}, m_{\mbox{{\scriptsize \emph{out}}}})=(1,1)$ [see
    Proposition~\ref{prop:c2} below for definitions]
    rather than  the non-generic $(1,2)$ points which have lower
    Hausdorff dimension~\cite{kn:FINR2}), all FBW paths through a
    small neighborhood of $(x,t)$ coalesce quickly with each
    other.
    \item The fact that, as already mentioned, a deterministic
    space-time point $(x,t)$ is the splice point of the unique FBW
    path passing through $(x,t)$ is quite natural in the
    stochastic flow context, as a consequence of Piterbarg's
    results. Things are contractive going forward in time because
    of the convergence of forward paths to {\em coalescing}
    Brownian motions. To see that things are expansive at any
    earlier time, $t_0<t$, along the flow line through $(x,t)$, we
    note that the flow contracts points $(x,t_0)$ by time $t$
    towards small regions which in the scaling limit converge to a
    locally finite set of locations which has zero probability of
    including the deterministic location $x$. Let $x_L$ and $x_R$
    respectively denote the two locations closest to $x$,
    respectively to its left and right, toward which points are
    contracted. Then, in the stochastic flow, in order to find the
    $(x_0,t_0)$ with $t_0<t$ which flows exactly into $(x,t)$ one
    must stay within the small expansive region separating the
    larger contractive regions which contract towards $x_L$
    and $x_R$.
    \item An interesting question raised by our results is whether
    the fact that {\em every} $\gamma$  in the FBW has a unique
    splice point
    $(x^*,t^*)$ where it changes from expansive behavior for (all)
    earlier times to contractive behavior for (most) future times
    has a natural interpretation for stochastic flows (on
    appropriate space and time scales) that goes beyond the
    discussion of the previous paragraph.
\end{enumerate}

 In the next two sections, we present
construction and characterization results for the FBW, and then in
Section~\ref{sec:conv}, a convergence result, which we apply to
the stochastic flows of homeomorphisms of~\cite{kn:H,kn:P} in
Section~\ref{sec:sf}. In Section~\ref{sec:max}, we discuss a
maximality property of the forward path of the FBW and its
relevance for the characterization of the (standard) Brownian web.
Finally, Section~\ref{app} contains proofs of two lemmas and a
proposition stated in earlier sections.

%%%%%%%%%%%%%%%%%%%%%%%%%%%%%%%%%%%%%%%%%%%%%%%%%%%%%%%%%%%%%%%%%%%%%%%%%%%%%%%
%%%%%%%%%%%%%%%%%% SECTION 2 %%%%%%%%%%%%%%%%%%%%%%%%%%%%%%%%%%%%%%%%%%%%%%%%%%%%%%%
%%%%%%%%%%%%%%%%%%%%%%%%%%%%%%%%%%%%%%%%%%%%%%%%%%%%%%%%%%%%%%%%%%%%%%%%%%%%%%%

\section{Construction}
\setcounter{equation}{0}
\label{sec:mbw}

In this section, we construct the FBW, combining backward and
forward paths of the double Brownian web (DBW).

%%%%%%%%%%%%%%%%%%%%%%%%%%%%%%%%%%%%%%%%%%%%%%%%%%%%%%%%%%%%%%%%%%%%%%%%%%%%%%%

\subsection{Preliminaries}
\label{ssec:pre}

As in~\cite{kn:FINR,kn:FINR1,kn:FINR2}, we begin with
$(\br^2,\rho)$, the completion (or compactification) of $\R^2$ under
the metric $\rho$, where
\begin{equation}
\label{rho}
\rho((x_1,t_1),(x_2,t_2))=
\left|\frac{\tanh(x_1)}{1+|t_1|}-\frac{\tanh(x_2)}{1+|t_2|}\right|
\vee|\tanh(t_1)-\tanh(t_2)|.
\end{equation}
$\br^2$ may be thought as
the image of $[-\infty,\infty]\times[-\infty,\infty]$
under the mapping
\begin{equation}
\label{compactify}
(x,t)\leadsto(\Phi(x,t),\Psi(t))
\equiv\left(\frac{\tanh(x)}{1+|t|},\tanh(t)\right).
\end{equation}

Let $\Pi^F$ denote the set of functions $f$ from 
$[-\infty,\infty]$ to $[-\infty,\infty]$ such that $\Phi(f(t),t)$
is continuous, and let
\begin{equation}
\label{d} d^F(f_1,f_2)= \sup_{-\infty\leq
t\leq\infty}|\Phi(f_1(t),t)-\Phi(f_2(t),t)|.
\end{equation}
be a metric on $\Pi^F$. Then $(\Pi^F,d^F)$ is a complete separable
metric space. We will throughout identify an arbitrary
element $f$ of $\Pi^F$ with its (space-time) path 
$\gamma=\{(f(s),s):\,-\infty\leq s\leq\infty\}$.

For the next result (whose proof we leave as an exercise), we use
the notation of Section 3 in~\cite{kn:FINR2}, where forward
semipaths $f$ starting at time $t_0$ are denoted $(f,t_0)$ and
belong to a space $\Pi$ with a similar notation for backward paths
$(g^b,t_0)$ in $\Pi^b$. Further, for
$[(f,t_0),(g^b,t_0)]\in\Pi\times\Pi^b$ such that
$f(t_0)=g^b(t_0)$, we let $(f,g^b,t_0)\in\Pi^F$ denote the path
through $(f(t_0),t_0)$ coinciding with $(f,t_0)$ after $t_0$, and
with $(g^b,t_0)$ before $t_0$.

\blem
   \label{lem:l1}
Let $\{[(f_n,t_n),(g^b_n,t_n)];\,n\geq1\}$ be a sequence of paths in
$\Pi\times\Pi^b$
%\footnote{For $n\geq1$, $(f_n,t_n)$ is a forward path starting
%in $(f_n(t_n),t_n)\in\bar\R^2$ and  $(g_n,t_n)$ is a forward path starting
%in $(g_n(t_n),t_n)\in\bar\R^2$.}
such that $f_n(t_n)=g^b_n(t_n)$ for $n\geq1$.
If $[(f_n,t_n),(g^b_n,t_n)]\to[(f,t),(g^b,t)]$ as $n\to\infty$ in
$\Pi\times\Pi^b$, then $(f_n,g^b_n,t_n)\to(f,g^b,t)$ as $n\to\infty$ in $\Pi^F$.
%Conversely, if $(f_n,g^b_n,t_n)\to(f,g^b,t)$ as $n\to\infty$ in $\Pi^F$
%and $(f_n(t_n),t_n)\to(f(t),t)$ in $\bar\R^2$, then
%$[(f_n,t_n),(g^b_n,t_n)]\to[(f,t),(g^b,t)]$ as $n\to\infty$ in
%$\Pi\times\Pi^b$.
\elem

Let now $\h^F$ denote the set of compact
subsets of $(\Pi^F,d^F)$, with $d_{\h^F}$ the induced Hausdorff metric, i.e.,
\begin{equation}
\label{dh}
d_{\h^F}(K_1,K_2)=\sup_{g_1\in K_1}\inf_{g_2\in K_2}d^F(g_1,g_2)\vee
        \sup_{g_2\in K_2}\inf_{g_1\in K_1}d^F(g_1,g_2).
\end{equation}
$(\h^F,d_{\h^F})$ is also a complete separable metric space.

The following lemmas will be used in Section~\ref{sec:char}; they
are key to our characterization and hence convergence results.
Their proofs are given in Section~\ref{app}.

\blem
   \label{lem:l2}
Let $\breve\W\in\h^F$ be noncrossing and such that for a dense countable
${\cal D}\subset\r^2$ and every $(x,t)\in{\cal D}$ there exists a unique path
$\gamma_{x,t}\in\breve\W$ passing through $(x,t)$. Then $\breve\W$ is determined
by the set of semipaths
$\W=\{W_{x,t}=(\gamma_{x,t}(s),s)_{s\geq t},\,(x,t)\in{\cal D}\}$.
\elem

Again we will use the notation of  Section 3 of~\cite{kn:FINR2}),
where $\h$ denotes the Hausdorff space of compact collections of
forward semipaths.

\blem
   \label{lem:l3}
For a noncrossing element $\breve\W\in\h^F$, let $\vec\W\in\h$ be the set of the forward
semipaths of $\breve\W$, i.e.,
$$\vec\W=\{W_{x,t}=(\gamma(s),s)_{s\geq t}:\,\gamma\in\breve\W,\,\gamma(t)=x,\,(x,t)\in\bar\R^2\}.$$
Suppose that for a dense countable
${\cal D}\subset\r^2$ and every $(x,t)\in{\cal D}$ there exists a unique semipath
$\gamma_{x,t}\in\vec\W$ starting at $(x,t)$, which furthermore avoids all
other points of ${\cal D}$. Then $\breve\W$ is determined
by the set of semipaths
$\W=\{W_{x,t},\,(x,t)\in{\cal D}\}$.
\elem

%\bigskip

%%%%%%%%%%%%%%%%%%%%%%%%%%%%%%%%%%%%%%%%%%%%%%%%%%%%%%%%%%%%%%%%%%%%%%%%%%%%%%%

\subsection{Two constructions}
\label{ssec:cons}

Let ${\cal D}$ be an (ordered) dense countable deterministic subset of $\r^2$.
We give two construction of the FBW, both based on the double Brownian web
(DBW).

\paragraph{\bf 1st Construction}
Let $\W^D=\{\tilde W_1,\tilde W_1^b,\tilde W_2,\tilde
W_2^b,\ldots\}$ be the DBW skeleton using the set ${\cal
D}=\{(x_1,t_1),(x_2,t_2),\ldots\}$ of space-time points as
starting points for these coalescing/reflecting forward/backward
Brownian motions (see the beginning of Section 3
of~\cite{kn:FINR2}). Define the FBW skeleton $\W^F=\{\tilde
W_1^F,\tilde W_2^F,\ldots\}$ by
\begin{equation}
   \label{eq:skel}
   \tilde W_j^F=
   \begin{cases}
     \tilde W_j^b(t),& t\leq t_j\\
     \tilde W_j(t),& t\geq t_j.
   \end{cases}
\end{equation}
Finally, let $\bw^F$ be the closure of $\W^F$ in $(\Pi^F,d^F)$.
\bprop
\label{prop:c1}
$\bw^F$ is compact (and so is an $(\h^F,d_{\h^F})$-valued random variable) and its
distribution does not depend on $\cal D$ or its ordering.
\eprop

\begin{figure}[!ht]
\begin{center}
\includegraphics[width=6cm]{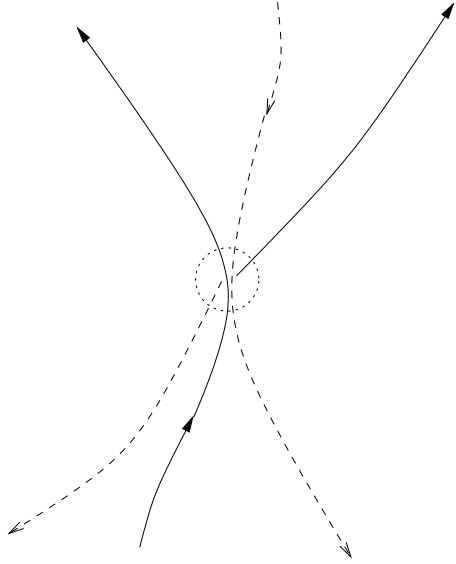}
\caption{A schematic diagram of a point $(x^*,t^*)$ of type
$(m_{\mbox{{\scriptsize\emph{in}}}},m_{\mbox{{\scriptsize
\emph{out}}}})=(1,2)$.}
\label{fig12}
\label{12} %don't know what this does
\end{center}
\end{figure}

\paragraph{\bf 2nd Construction}
As in~\cite{kn:FINR2}, we denote by $\bar\W$ and $\bar\W^b$ the
forward and backward Brownian webs, which are compact subsets of
semipaths of their respective spaces $\Pi$ and $\Pi^b$. \bprop
\label{prop:c2} $\bw^F$ consists of the two ``trivial'' paths that
are identically $+\infty$ or $-\infty$ for all $t$, plus the
following collection of paths constructed from the DBW
$\bw^D=(\bw,\bw^b)$. For each $(x^*,t^*)\in\r^2$, of type
$(m_{\mbox{{\scriptsize\emph{in}}}},
m_{\mbox{{\scriptsize\emph{out}}}})$\footnotemark, let
$W_1,\ldots,W_{\!\!m_{\mbox{{\tiny\emph{out}}}}}$ be the forward
paths from $(x^*,t^*)$, and
$W_1^b,\ldots,W^b_{\!\!m_{\mbox{{\tiny\emph{in}}}}+1}$ be the
backward paths from $(x^*,t^*)$. Then take every path of the form
\begin{equation}
   \label{eq:splice}
   W_{ij}^F=(W_i,W_j^b,t^*)=
   \begin{cases}
     W_i(t),& t\geq t^*\\
     W_j^b(t),& t\leq t^*,
   \end{cases}
\end{equation}
{\em except}, when
$(m_{\mbox{{\scriptsize\emph{in}}}},m_{\mbox{{\scriptsize\emph{out}}}})=(1,2)$
(see Figure~\ref{fig12}), do not take the {\em unique} choice of
$i,j$ for which $W_{ij}^F$ would cross the forward path passing
{\em through} $(x^*,t^*)$ [i.e., in Figure~\ref{fig12}, do not
take the $W_{ij}^F$ from Southwest to Northeast]. \eprop

\footnotetext{$m_{\mbox{{\scriptsize\emph{in}}}}$ denotes the
number of forward paths of $\bw^D$ starting before $t^*$, not
coalescing up to right before $t^*$, and passing through
$(x^*,t^*)$, and $m_{\mbox{{\scriptsize\emph{out}}}}$ denotes the
number of forward paths of $\bw^D$ starting at $(x^*,t^*)$; see,
e.g., Section 3 of~\cite{kn:FINR2}.}

%%%%%%%%%%%%%%%%%%%%%%%%%%%%%%%%%%%%%%%%%%%%%%%%
%%%%%%%%%%%%%%%%%%  FIGURE  %%%%%%%%%%%%%%%%%%%%%%%%
%%%%%%%%%%%%%%%%%%%%%%%%%%%%%%%%%%%%%%%%%%%%%%%%

%\eject

\brm \label{rmk:trivial} The trivial paths, namely the identically
$-\infty$ one and the identically $+\infty$ one, are also
splicings of paths in the DBW. Indeed, for any $t^*$, the
identically $+\infty$ path is the splicing of the forward path (in
$\bar\W$) starting at $(+\infty,t^*)$ with the backward path (in
$\bar\W^b$) from $(+\infty,t^*)$. One difference between
the trivial and the nontrivial paths is that the latter have a
well-defined unique splice point (in $\R^2$), while for the latter
all their points can be seen as splice points. \erm

\brm \label{rmk:splice} Every path $\g$ in the FBW consists of a
backward path from the DBW between times $t=-\infty$ and $t=$ some
$t^*(\g)$, spliced with a forward path from $t^*(\g)$ to
$t=+\infty$. {\em Every} $(x^*,t^*)$ in $\r^2$ is a splice point
of of one or more $\g$'s from the FBW; the number of $\g$'s for a
given splice point $(x^*,t^*)$, determined by the type
$(m^*_{\mbox{{\scriptsize\emph{in}}}},m^*_{\mbox{{\scriptsize\emph{out}}}})$
of $(x^*,t^*)$, is
$(m^*_{\mbox{{\scriptsize\emph{in}}}}+1)\,m^*_{\mbox{{\scriptsize\emph{out}}}}$,
except that for a point of type $(1,2)$, that number is only $3$,
rather than $4$.

As one follows a backward path $\g^b$ {\em forward} in time to
some $(x^\#,t^\#)$, the number of choices  of how to continue it
(locally) as a path in the FBW depends on the type
$(m^\#_{\mbox{{\scriptsize\emph{in}}}},m^\#_{\mbox{{\scriptsize\emph{out}}}})$
of $(x^\#,t^\#)$, as well as on whether $\g^b$ is part of a
backward path {\em passing through} $(x^\#,t^\#)$ [rather than
originating at $(x^\#,t^\#)$]. If $\g^b$ is not part of a backward
path passing through $(x^\#,t^\#)$ (this is the case if
$m^\#_{\mbox{{\scriptsize\emph{out}}}}=1$, or sometimes when
$(m^\#_{\mbox{{\scriptsize\emph{in}}}},m^\#_{\mbox{{\scriptsize\emph{out}}}})=(1,2)$
--- e.g.~if $\g^b$ is the Southwest path in Figure~\ref{fig12},
then $(x^\#,t^\#)$ must be a splice point and there is a unique
choice of forward path to splice into.  If $\g^b$ is rather part
of a backward path passing through $(x^\#,t^\#)$, then
$(x^\#,t^\#)$ can be a splice point with
$m^\#_{\mbox{{\scriptsize\emph{out}}}}$ choices of forward path to
splice into, or $\gamma^b$ can be continued to any of
$m^\#_{\mbox{{\scriptsize\emph{out}}}}-1$ possible backward path
continuations. E.g., for
$(m^\#_{\mbox{{\scriptsize\emph{in}}}},m^\#_{\mbox{{\scriptsize\emph{out}}}})=(0,3)$,
there are $5$ possible continuations, with $3$ of them
corresponding to splice points. Once a splice point has been
chosen, and a continuation to a forward path has been made, that
forward path is followed until $t=+\infty$ with no further
choices.

Note that if a path $\g$ in the FBW touches a point
$(\g(t_0),t_0)$ with $m_{\mbox{{\scriptsize\emph{out}}}}=1$, then
for $t\geq t_0$, $\g(t)$ follows the unique forward path in the BW
from $(\g(t_0),t_0)$. \erm

%\brm
%\label{rmk:det}
%Lemma~\ref{lem:l3} implies that for any deterministic dense countable
%${\cal D}\subset\R^2$,
%$\bar\W^F$ (playing the role of  $\breve\W$ in the lemma), is almost
%surely determined by the set of its semipaths starting from ${\cal D}$
%(corresponding to $\W$ in the lemma), which in this case is a
%Brownian web skeleton, and is readily seen to satisfy the conditions
%of the lemma.
%\erm

\noindent{\bf Proof of Proposition~\ref{prop:c2}} Let $\hat\W$ be
the set of paths described in the statement. Also let $\bar\W^F$
denote the closed set of paths in $(\Pi^F,d^F)$ from our first
construction for some particular $\cd$; as we have not 
yet proved
Proposition~\ref{prop:c1}, $\bar\W^F$ could a priori depend on
$\cd$. We show $(i)\,\bar\W^F\subset\hat\W$ and
$(ii)\,\hat\W\subset\bar\W^F$.

$(i)$ Take a path $W^F\in\bar\W^F$ and a sequence $\tilde
W_{j_1}^F,\tilde W_{j_2}^F,\ldots$ of paths in $\W^F$ converging
to $W^F$. By the compactness of $\bar\R^2$, the sequence of splice
points $(x_{j_1},t_{j_1}),(x_{j_2},t_{j_2}),\ldots$ of $\tilde
W_{j_1}^F,\tilde W_{j_2}^F,\ldots$ has a convergent subsequence
$(x_{j'_1},t_{j'_1}),(x_{j'_2},t_{j'_2}),\ldots$ with a limit
point $(x^*,t^*)\in\bar\R^2$. It readily follows that $(\tilde
W_{j'_i})_{i\geq1}$ converges in $\Pi$ to a forward path $W$
starting in $(x^*,t^*)$, that $(\tilde W^b_{j'_i})_{i\geq1}$
converges in $\Pi^b$ to a backward path $W^b$ starting in
$(x^*,t^*)$, and that $W^F$ is the splicing of $W$ and $W^b$. By
the compactness of the double Brownian web, $W\in\bar\W$ and
$W^b\in\bar\W^b$. It is also easy to see that $(W,W^b,t^*)$ cannot
be the one case of splicing ruled out in the definition of
$\hat\W$. Thus, $W^F\in\hat\W$. (We note that $W^F$ is a trivial
path iff $(x^*,t^*)\in\bar\R^2\setminus\R^2$.)

$(ii)$ Take a path $\hat W^F\in\hat\W$, and let $(x^*,t^*)$ be its
splice point (taken to be $(\pm\infty,0)$ for the identically
$\pm\infty$ trivial paths, respectively), so that $\hat W^F=(\hat
W,\hat W^b,t^*)$.

If $(x^*,t^*)$ is {\em not} of type $(1,2)$\footnote{This happens
in particular in the case of $(\pm\infty,0)$, which is of type
$(1,1)$.}, then either $\hat W$ is the unique path in $\bar\W$
from $(x^*,t^*)$, or $\hat W^b$ is the unique path in
$\bar\W^b$ starting from $(x^*,t^*)$. Suppose we are in the the
latter case (the former case can be handled similarly). Then,
since $\bar\W$ is the closure of the set of its paths starting
from points in $\cd$, there exists a sequence
$(x_{j_1},t_{j_1}),(x_{j_2},t_{j_2}),\ldots$ in ${\cal D}$
approaching $(x^*,t^*)$ such that $(\tilde W_{j_i})_{i\geq1}$
converges in $\Pi$ to $\hat W$. Now by the compactness of
$\bar\W^b$ and the uniqueness property of $\hat W^b$ just
mentioned, $(\tilde W^b_{j_i})_{i\geq1}$ converges in $\Pi^b$ to
$\hat W^b$. By Lemma~\ref{lem:l1}, $(\tilde W^F_{j_i})_{i\geq1}$
converges in $\Pi^F$ to $\hat W^F$, and thus $\hat
W^F\in\bar\W^F$.

Suppose now that $(x^*,t^*)$ is a left-handed type $(1,2)$ point
as in Figure~\ref{fig12} (see Remark 3.11 in~\cite{kn:FINR2}; the
right-handed case is treated similarly), and let $W_1$, $W_2$ be
the two paths in $\bar\W$ starting from $(x^*,t^*)$, and $W_1^b$,
$W_2^b$ be the two paths in $\bar\W^b$ starting from $(x^*,t^*)$,
with $W_1$ to the left of $W_2$, and $W_1^b$ to the left of
$W_2^b$. Let also $\check W_1$ be any path in $\bar\W$ starting
before $t^*$ and passing through $(x^*,t^*)$, and $\check W_2^b$
be any path in $\bar\W^b$ from above $t^*$ passing through
$(x^*,t^*)$. ($\check W_1$ and $\check W_2^b$ coincide with $W_1$
and $W_2^b$ above and below $t^*$ respectively.) We then have that
$(\hat W,\hat W^b,t^*)=(W_i,W_j^b,t^*)$, for some $i,j=1,2$, but
$(i,j)\ne(2,1)$.

If $(i,j)=(1,1)$ and we take a sequence
$(x_{j_1},t_{j_1}),(x_{j_2},t_{j_2}),\ldots$ in ${\cal D}$ approaching
$(x^*,t^*)$ such that $(x_{j_i},t_{j_i})$ is to the left of $\check W_1$
for all $i$, then $\tilde W_{j_i}\to W_1$, and $\tilde W^b_{j_i}\to W^b_1$
as $i\to\infty$.

If $(i,j)=(2,2)$ and we take a sequence
$(x_{j'_1},t_{j'_1}),(x_{j'_2},t_{j'_2}),\ldots$ in ${\cal D}$
approaching $(x^*,t^*)$ such that $(x_{j'_i},t_{j'_i})$ is to the
right of $\check W_2^b$ for all $i$, then $\tilde W_{j'_i}\to
W_2$, and $\tilde W^b_{j'_i}\to W^b_2$ as $i\to\infty$.

If $(i,j)=(1,2)$, then if we take a sequence
$(x_{j''_1},t_{j''_1}),(x_{j''_2},t_{j''_2}),\ldots$ in ${\cal D}$ approaching
$(x^*,t^*)$ such that $(x_{j''_i},t_{j''_i})$ is between $\check W_1$ and
$\check W_2^b$ for all $i$, then $\tilde W_{j''_i}\to W_1$, and
$\tilde W^b_{j''_i}\to W^b_2$ as $i\to\infty$.

In all cases, we conclude from Lemma~\ref{lem:l1} that $\hat W^F$ is the limit
of a sequence in $\bar\W^F$. $\square$

\bigskip

\noindent{\bf Proof of Proposition~\ref{prop:c1}} Let $W^F_n=(W_n,W^b_n,t_n)$,
$n\geq1$, be a sequence of paths in $\bar\W^F$ with splice points
$(x_n,t_n)$, $n\geq1$, respectively. By the compactness of $\bar\R^2$,
$(x_{n'},t_{n'})\to(x,t)\in\bar\R^2$ as $n'\to\infty$ for a subsequence
$(x_{n'},t_{n'})$ of $(x_{n},t_{n})$. By the compactness of $\bar\W^D$,
$(W_{n''},W^b_{n''})\to(W,W^b)$ as $n''\to\infty$
for a a subsequence $(x_{n''},t_{n''})$ of $(x_{n'},t_{n'})$, with
$(W,W^b)\in\bar\W^D$. Since clearly $W(t)=W^b(t)=x$, we have by
Lemma~\ref{lem:l1} that $W^F_{n''}\to(W,W^b,t)\in\bar\W^F$ as $n''\to\infty$.
The nondependence of the distribution of $\bar\W^F$ on ${\cal D}$ is immediate
from the same property of $\bar\W^D$ and Proposition~\ref{prop:c2}.
$\square$

%%%%%%%%%%%%%%%%%%%%%%%%%%%%%%%%%%%%%%%%%%%%%%%%%%%%%%%%%%%%%%%%%%%%%%%%%%%%
%%%%%%%%%%%%%%%% SECTION 3 %%%%%%%%%%%%%%%%%%%%%%%%%%%%%%%%%%%%%%%%%%%%%%%%%%%%%%
%%%%%%%%%%%%%%%%%%%%%%%%%%%%%%%%%%%%%%%%%%%%%%%%%%%%%%%%%%%%%%%%%%%%%%%%%%%%%

\section{Characterization}

\setcounter{equation}{0}
\label{sec:char}

In this section we state some results characterizing the
distribution of the FBW.  We will use the following fact about the
double Brownian web, whose proof is found in Section~\ref{app}.

\bprop \label{prop:dbw} Let $\bar\W\in(\h,\f_{\h})$ be a standard
Brownian web. Then there exists a standard double Brownian web
$\hat\W^D=(\hat\W,\hat\W^b)\in(\h^D,\f_{\h^D})$ such that
a.s.~$\hat\W=\bar\W$. \eprop

\bteo
 \label{teo:char1}
There is an $(\h^F,d_{\h^F})$-valued random variable $\bw^F$ whose distribution is
uniquely determined by the following properties.

($a$) Almost surely the paths of $\bw^F$ are noncrossing (although
they may touch, including coalescing and bifurcating).

($b_1$) From any deterministic point $(x,t)\in\r^2$, there is
almost surely a unique path $W_{x,t}^F$ passing through $x$ at
time $t$.

($b_2$) For any deterministic $n$, $\{(x_1,t_1),\ldots,(x_n,t_n)\}$, the joint distribution
of the semipaths $\{W_{x_j,t_j}^F(t),\,t\geq t_j,\,j=1,\ldots,n\}$ is that of
coalescing Brownian motions (with unit diffusion constant).
\eteo

\brm
\label{rmk:char1}
($b_1$) and ($b_2$) may together be replaced by the following weaker condition.

(b) For a given deterministic countable dense set ${\cal D}\in\r^2$, there exist
(not necessarily unique a priori) paths $W_{x_j,t_j}^F$ in $\bw^F$ for each
$(x_j,t_j)$ in ${\cal D}$ such that
$\{W_{x_j,t_j}^F(t),\,t\geq t_j,\,j=1,2,\ldots\}$ is distributed as coalescing
Brownian motions (with unit diffusion constant) starting at the points of $\cal D$.
\erm

\noindent{\bf Proof of Theorem~\ref{teo:char1}} Condition ($b_2$)
implies that for any fixed deterministic countable dense ${\cal
D}$, the set of semipaths $\W=\{(W_{x,t}^F(s))_{s\geq
t},\,(x,t)\in{\cal D}\}$ is distributed as a Brownian web
skeleton. Thus $\bar\W$, the closure of $\W$, is a standard
Brownian web. By Proposition~\ref{prop:dbw} there exists a dual
Brownian web $\bar\W^b$ such that $(\bar\W,\bar\W^b)$ is a double
Brownian web. Condition ($b_1$) implies that there is a unique
backward semipath of $\bar\W^F$ starting from each point of ${\cal
D}$. Thence Lemma~\ref{lem:l2} and Proposition~\ref{prop:c2} (or
Proposition~\ref{prop:c1}) yield the result. $\square$

\bigskip

\noindent{\bf Proof of Remark~\ref{rmk:char1}} Condition ($b$)
implies that the set of semipaths $\W=\{(W_{x,t}^F(s))_{s\geq
t},\,(x,t)\in{\cal D}\}$ is a Brownian web skeleton (the
a.s.~uniqueness of the semipath starting in each point of ${\cal
D}$ follows by the {\em trapping argument} used in the proof of
Proposition 3.1 of~\cite{kn:FINR1} (arXiv version)). The proof
continues like the one of Theorem~\ref{teo:char1} above, using
Lemma~\ref{lem:l3} instead of Lemma~\ref{lem:l2}. $\square$

%\bteo
% \label{teo:char2} There is an $(\h^F,d_{\h^F})$-valued random
%variable $\bw^F$ whose distribution is uniquely determined by (b)
%and the following minimality property.

%($a'$) If $\W^F_*$ is any other $(\h^F,d_{\h^F})$-valued random
%variable satisfying (b), then $\mu_{\bw^F}<<\mu_{\W^F_*}$, where
%$\mu_{\bw^F},\,\mu_{\W^F_*}$ are the probability distributions of
%$\bw^F,\,\W^F_*$, respectively, and $<<$ denotes the stochastic
%ordering corresponding to the $\subseteq$ partial ordering of
%$\h^F$. \eteo

%\noindent{\bf Proof} Following the arguments in the proof of
%Theorem~\ref{teo:char2} above, $\bw^F$ contains the FBW, denoted
%$\hat\W^F$, obtained from $(\bar\W,\bar\W^b)$ as in Construction 2
%(see proof of Theorem~\ref{teo:char1} above). Since $\hat\W^F$
%also satisfies $(b)$, $(a')$ implies that
%a.s.~$\hat\W^F=\bar\W^F$. $\square$

%%%%%%%%%%%%%%%%%%%%%%%%%%%%%%%%%%%%%%%%%%%%%%%%%%%%%%%%%%%%%%%%%%%%%%%%%%%%
%%%%%%%%%%%%%%%% SECTION 4 %%%%%%%%%%%%%%%%%%%%%%%%%%%%%%%%%%%%%%%%%%%%%%%%%%%%%%
%%%%%%%%%%%%%%%%%%%%%%%%%%%%%%%%%%%%%%%%%%%%%%%%%%%%%%%%%%%%%%%%%%%%%%%%%%%%%

\section{Convergence}

\setcounter{equation}{0}
\label{sec:conv}

In this section we establish convergence criteria for a sequence
of $(\h^F,d_{\h^F})$-valued random variables to converge to the
FBW. These criteria are then applied in the next section to show
that the diffusively rescaled stochastic flow of homeomorphisms
of~\cite{kn:H,kn:P} converges to the FBW as the scale parameter
goes to $0$.

\bteo
\label{teo:conv}
Let $\X_1,\X_2,\ldots$ be a sequence of $(\h^F\!,d_{\h^F})$-valued random variables
such that

($\bar a$) almost surely, the paths of each $\X_m$ are noncrossing, and

($\bar b$) for a given deterministic countable dense set ${\cal
D}\in\r^2$, there exists for each $(x_j,t_j)\in{\cal D}$ a path
$\theta^{x_j,t_j}_m$ in $\X_m$ such that the collection of
semipaths $\V_m:=\{(\theta^{x_j,t_j}_m(t),\,t\geq t_j),
(x_j,t_j)\in{\cal D}\}$ converges in distribution as $m\to\infty$
to coalescing Brownian motions (with unit diffusion constant)
starting in ${\cal D}$.

Then $\X_m$ converges in distribution as $m\to\infty$ to the FBW.

\eteo

\noindent{\bf Proof}\,  We claim that $(\X_m,\,m\geq1)$ is tight.
Then, conditions ($\bar a)$-$(\bar b$) and Remark~\ref{rmk:char1}
imply that all limit points of $(\X_m,\,m\geq1)$ are distributed
as the FBW.

To justify the tightness claim, let $\V_m$ (resp., $\V^b_m$)
denote the collection of all forward (resp., backward) semipaths
of $\X_m$. By Proposition B.3 of~\cite{kn:FINR1}, we have
tightness of $\V_m$ in $(\h,d_\h)$ and it then suffices to verify
that $(\V^b_m,\,m\geq1)$ is tight in $(\h^b,d_{\h^b})$. This
readily implies tightness of $\{(\V_m,\V^b_m),\,m\geq1\}$ in
$(\h^D,d_{\h^D})$, and this and Lemma~\ref{lem:l1} imply tightness
of $(\X_m,\,m\geq1)$.

Now tightness of $(\V^b_m,\,m\geq1)$ follows by the very same {\em
blocking argument} used in the proof of Proposition~\ref{prop:dbw}
given in Section~\ref{app} below to verify the tightness condition
$(T_1)$ from Appendix B of~\cite{kn:FINR1}. $\square$

%\eject

%%%%%%%%%%%%%%%%%%%%%%%%%%%%%%%%%%%%%%%%%%%%%%%%%%%%%%%%%%%%%%%%%%%%%%%%%%%%
%%%%%%%%%%%%%%%% SECTION 5 %%%%%%%%%%%%%%%%%%%%%%%%%%%%%%%%%%%%%%%%%%%%%%%%%%%%%%
%%%%%%%%%%%%%%%%%%%%%%%%%%%%%%%%%%%%%%%%%%%%%%%%%%%%%%%%%%%%%%%%%%%%%%%%%%%%%

\section{Stochastic flows and the FBW}

\setcounter{equation}{0}
\label{sec:sf}

Let $\Xi=\{\xi_{st},\,s\leq t\}$ be an isotropic stochastic flow
of homeomorphisms with covariance structure $B$. This means, among
other things, that for every $s\in\R$
$$(\xi_{st}(x_1),\ldots,\xi_{st}(x_n);\,t\geq s)$$
is an $R^n$-diffusion starting at $(x_1,\ldots,x_n)$
with generator
$$\frac{B(0)}{2}\sum^n\frac{\partial^2}{\partial
y_i^2}+\sum^n_{i<j}B(y_i-y_j)\frac{\partial^2}{\partial
y_i\partial y_j}.$$

{Piterbarg}~\cite{kn:P1,kn:P} studied expansions/contractions of
$\Xi$ via the $\d\to0$ limit of
$$\xi^\d_t(x)\equiv\d\,\xi_{0\,
t\d^{-2}}(x\d^{-1}).$$ It follows from his results that when $B$ is nice
(with $B(x)\to0$ as $|x|\to\infty$; see conditions (B1-B4) in~\cite{kn:P}),
then for all $(x_1,t_1),\ldots,(x_n,t_n)\in\R^2$
\beq
\label{pit}
\{(\xi^\d(x_1,t_1),\,s\geq t_1),\ldots,(\xi^\d(x_n,t_n),\,s\geq t_n)\}
\eeq
converges weakly as $\d\to0$ to coalescing Brownian motions starting from
$(x_1,t_1),\ldots,(x_n,t_n)\in\R^2$, where for $i=1,\ldots,n$ and $\d>0$,
$$\xi^\d(x_1,t_1)=\xi_{t_i\d^{-2}, s\d^{-2}}(x_i\d^{-1})$$
(see Theorem 6 in~\cite{kn:P}; there the case $t_i=$ constant is stated,
from which the general case readily follows).

Since these flows are noncrossing, the above and
Theorem~\ref{teo:conv} immediately imply the following result.
\bteo \label{teo:sf} For $\d>0$, let
$\X_\d=\{[\d\,\xi_{0\,t\d^{-2}}(x\d^{-1}),\,t\in\R],\,x\in\R\}$ as
in~\cite{kn:P}. Then $\X_\d$ converges in distribution to the FBW
as $\d\to0$. \eteo

%%%%%%%%%%%%%%%%%%%%%%%%%%%%%%%%%%%%%%%%%%%%%%%%%%%%%%%%%%%%%%%%%%%%%%%%%%%%
%%%%%%%%%%%%%%%% SECTION 6 %%%%%%%%%%%%%%%%%%%%%%%%%%%%%%%%%%%%%%%%%%%%%%%%%%%%%%
%%%%%%%%%%%%%%%%%%%%%%%%%%%%%%%%%%%%%%%%%%%%%%%%%%%%%%%%%%%%%%%%%%%%%%%%%%%%%

\section{The forward full Brownian web}

\setcounter{equation}{0}
\label{sec:max}

In this section, we use the FBW, or rather its forward part, to
provide an example clarifying the significance of a condition in 
characterization results for the (standard) Brownian web (see,
e.g., Condition ($ii$) in Theorem~2.1 or Condition ($ii'$)
in Theorem~4.1 [Theorem~4.5 in the arXiv version]
of~\cite{kn:FINR1} and the
remark following Theorem~2.1 there; see also Remark~\ref{rmk:cond} below).

For a path $\gamma\in\W^F$ and $t\in[-\infty,\infty]$, let
$\gamma_t$ denote the forward semipath starting at
$(\gamma(t),t)$, and consider the path collection
\begin{equation}
\label{max}
  \W^{FF}:=\{\gamma_t:\,\gamma\in\W^F,\,t\in[-\infty,\infty]\}.
\end{equation}
The next theorem says that $\W^{FF}$ is, in a natural sense, the
{\em maximal} collection of (noncrossing) forward paths containing the
standard Brownian web; thus it is the {\em full forward} as well
as the {\em forward full} Brownian web.
Indeed, $\W^{FF}$ is the maximal collection
satisfying properties (0)-(iii) below while the standard BW is the
{\em minimal} such collection. In particular, for 
$-\infty < s < t < \infty$, the set of $x$ such that there is a
(forward) semipath passing through both $\r \times \{s\}$
and $(x,t)$
is locally finite for the BW while it consists of all of $\r$
for $\W^{FF}$.

\bteo \label{teo:max} $\W^{FF}$ has the following properties.
\begin{itemize}
\item[(o)] $\W^{FF}$ is a.s.~compact in $(\Pi,d)$,
           the space of forward semipaths\footnotemark;
\item[(i)] the paths of $\W^{FF}$ are noncrossing; \item[(ii)]
from any deterministic point \((x,t)\) in $\r^{2}$,
           there is almost surely a unique path \({W}_{x,t}\) in $\W^{FF}$
           starting from \( (x,t) \);
\item[(iii)] for any deterministic \( n, (x_{1}, t_{1}), \ldots,
           (x_{n}, t_{n}) \), the joint distribution of \(
           {W}_{x_{1},t_{1}}, \ldots, {W}_{x_{n},t_{n}} \) is that
           of coalescing Brownian motions (with unit diffusion constant);
\item[(iv)] $\W^{FF}$ strictly contains a version of the standard
Brownian web $\bar\W$; \item[(v)] if $\hat\W$ satisfies properties
(o)-(iii) above, then there exists
           a version $\W^\ast$ of $\W^{FF}$ such that $\hat\W\subset\W^\ast$.
\end{itemize}
\eteo

\footnotetext{See Section 3 of~\cite{kn:FINR1} for the precise
definition of the metric $d$ for semipaths.}

\brm\label{rmk:cond} Theorem 2.1 of~\cite{kn:FINR1} establishes
that under conditions (ii) and (iii) of Theorem~\ref{teo:max}
above, a $(\h,d_\h)$-random variable is a (standard) Brownian web
provided that a further condition
(alluded to at the beginning of this section)
is satisfied, namely that it
is almost surely the closure of the set of its paths starting at a
deterministic countable subset of $\bar\r^2$.
Theorem~\ref{teo:max} above makes clear that without the
latter condition, we can have $(\h,d_\h)$-random variables
satisfying conditions (ii) and (iii) of Theorem~\ref{teo:max}
other than the Brownian web (as established in 
Theorem~3.1 [Theorem~3.6 in the arXiv version] 
of~\cite{kn:FINR1}, these $(\h,d_\h)$-random variables must be
stochastically bigger than the Brownian web; that is of course the
case with $\W^{FF}$, as stated in (iv) above).\erm

\noindent{\bf Proof} Properties (o)-(iii) have already been
discussed in/follow readily from our construction (e.g., the 2nd
one).

(iv) That $\W^{FF} \supset \bar\W$ is immediate
from our 2nd construction. That $\W^{FF}$ has 
strictly more paths than $\bar\W$ is seen, e.g., in
that there are infinitely many paths of $\W^{FF}$ starting at any
$(0,2)$-point of $\bar\W$, namely, apart from the 
two paths $W_1$ and
$W_2$, say, of $\bar\W$ starting at such a point, all forward
semipaths of $\W^F$ obtained by splicing paths of $\bar\W$
starting at points above the $(0,2)$-point 
and between $W_1$ and
$W_2$, with the path segments in the dual to $\bar\W$ starting at
those points and going down to the $(0,2)$-point.

(v) For a deterministic dense countable $\cd\subset\R^2$, let $\W$
be the set of paths of $\hat\W$ starting from $\cd$. Then, by
(ii)-(iii), $\bar\W$, the closure of $\W$, is a Brownian web. By
Proposition~\ref{prop:dbw} there exists a dual Brownian web
$\bar\W^b$ such that $\bar\W^D=(\bar\W,\bar\W^b)$ is a DBW. Let
now $\bar\W^F$ be the FBW obtained from $\bar\W^D$ (as in our 2nd
construction), and $\bar\W^{FF}$ be the 
forward full Brownian web
obtained from $\bar\W^F$ (as in~(\ref{max}) above). We claim that
$\hat\W\subset \bar\W^{FF}$. To justify that, we will {\em extend}
$\hat\W$ to a subset of $\h^F$ and then use (o)-(i) and
Lemma~\ref{lem:l2}.

To make that extension, for each $(x,t)\in \r^2$
we choose a single backward path $\hat
W^b_{x,t}$ in $\bar\W^b$ from $(x,t)$, with the
proviso that if $(x,t)$ is a $(1,2)$-point of $\bar\W^D$, then
this single backward path must be chosen in the unique way so that
the splicings of $\hat W^b_{x,t}$ with the forward paths of
$\bar\W$ from $(x,t)$ do not cross 
any of the paths of $\bar\W^D$
passing through $(x,t)$ (e.g., for the $(1,2)$-point of
Figure~\ref{fig12}, we choose the Southwest backward path). For
each path 
$\hat W\in\hat\W$
starting at $(x,t)$, let $\g=\g(\hat W)$ be the path
in $\Pi^F$ obtained as the splicing of $\hat W$ and $\hat
W^b_{x,t}$, and consider the sets 
$\W^*=\{\g(\hat W):\,\hat W\in\hat\W\}$,
and its closure $\bar\W^*$. We claim that $\W^*$ is noncrossing,
since $\hat\W$ is noncrossing by (i), and if $\hat W\in\hat\W$
crossed any $\hat W^b_{x,t}\in\bar\W^b$, then there would have to
be a path in $\W$ crossing $\hat W^b_{x,t}$, which cannot happen
since $\W$ is part of the DBW $\bar\W^D$; we have also to
rule out that $\hat W\in\hat\W$ crosses a path
$\g\in\W^*$ at a splicing point $(x,t)$ of $\g$, but that would
mean that $(x,t)$ is a $(1,2)$-point of $\bar\W^D$ and that the
choice of $\hat W^b_{x,t}$ did not obey the above proviso.
Compactness of $\bar\W^*$ follows readily from (o) and the
precompactness of $\{\hat W^b_{x,t},\,(x,t)\in\bar\R^2\}$. The
claim of (v) readily follows. $\square$

%\footnotetext{See Section 3 of~\cite{kn:FINR1}.}

%%%%%%%%%%%%%%%%%%%%%%%%%%%%%%%%%%%%%%%%%%%%%%%%%%%%%%%%%%%%%%%%%%%%%%%%%%%%
%%%%%%%%%%%%%%%% APPENDIX %%%%%%%%%%%%%%%%%%%%%%%%%%%%%%%%%%%%%%%%%%%%%%%%%%%%%%
%%%%%%%%%%%%%%%%%%%%%%%%%%%%%%%%%%%%%%%%%%%%%%%%%%%%%%%%%%%%%%%%%%%%%%%%%%%%%

%\appendix

\section{Several proofs}
\setcounter{equation}{0}
\label{app}

We begin with a Lemma that will be used in the proof of
Lemma~\ref{lem:l2}.
\blem
\label{lm:close}
Under the hypotheses of
Lemma~\ref{lem:l2}, $\breve\W$ is the closure of
$\{\gamma_{x,t}:\,(x,t)\in\cd\}$.
\elem

\noindent{\bf Proof} Let $\G$ denote
$\{\gamma_{x,t}:\,(x,t)\in\cd\}$ and take any path $\g\in\breve\W$.
The lemma follows if we can show that $\gamma$ is in the closure
of $\Gamma$.
Let now $\{q_i,\,i=1,2,\ldots\}$ be an enumeration of the rational
numbers of $\R$. We start by taking a sequence 
$(s^{(1)}_i)_{i=1}^\infty$ 
such that for $i=1,2,\ldots$ we have $s^{(1)}_i=(x^{(1)}_i,t^{(1)}_i)\in\cd$, 
$x^{(1)}_i\geq\gamma(t^{(1)}_i)$ and $s^{(1)}_i\to(\gamma(q_1),q_1)$ as $i\to\infty$.
Let $\g^{(1)}$ be a limit (in $\Pi^F$) of $\g_{x^{(1)}_i,t^{(1)}_i}$.
It is clear that $\g^{(1)}$ is in $\breve\W$, that
$\g^{(1)}(s) \geq \g(s) \,\, \forall \, s \in \r$, and that 
$\g^{(1)}(q_1)=\gamma(q_1)$.

Let us suppose that for $n\geq1$, we have chosen paths $\g^{(j)}$, 
$j=1,\ldots,n$, such that  
\begin{equation}
\label{eq:cond}
\g^{(j)}\in\breve\W,\,\g^{(j)}(s)\geq\gamma(s)\,\forall s\in\R,\,\mbox{ and }\,
\g^{(j)}(q_i)=\gamma(q_i)\,\mbox{ for }\,1\leq i\leq j\leq n;
\end{equation}
we now obtain a path $\g^{(n+1)}$ such~(\ref{eq:cond}) holds with
$n$ replaced by $n+1$. If $\g^{(n)}(q_{n+1})=\gamma(q_{n+1})$, then
let $\g^{(n+1)}=\g^{(n)}$; otherwise, we must have  
$\g^{(n)}(q_{n+1})>\gamma(q_{n+1})$. We can then choose a sequence
$s^{(n+1)}_i=(x^{(n+1)}_i,t^{(n+1)}_i)$ in $\cd$ such that
$\gamma(t^{(n+1)}_i)\leq x^{(n+1)}_i<\gamma^{(n)}(t^{(n+1)}_i)$ and 
$s^{(n+1)}_i\to(\gamma(q_{n+1}),q_{n+1})$ as $i\to\infty$. Now let
$\g^{(n+1)}$ be a limit (in $\Pi^F$) of $\g_{x^{(n+1)}_i,t^{(n+1)}_i}$.

Finally, let $\bar\gamma$ be a limit 
(in $\Pi^F$) of $\{\g^{(n)}:\,n\geq1\}$.
Then we must have $\bar\gamma(s) \geq \gamma(s) \,\,
\forall s \in \r$ while $\bar\gamma(q_i)=\gamma(q_i)$ for $i=1,2,\ldots$,
and thus $\bar\gamma=\gamma$.
$\square$

\bigskip

\noindent{\bf Proof of Lemma~\ref{lem:l2}} Write ${\cal
D}=\{(x_i,t_i),\,i=1,2,\ldots\}$ and $W_i=W_{x_i,t_i}$,
$i=1,2,\ldots$. By Lemma~\ref{lm:close}, it is enough to argue
that for every $i=1,2,\ldots$ the backward semipath
$W^b_i=(\gamma_{x_i,t_i}(s),s)_{s\leq t_i}$ is determined by
$\W=\{W_i,\,i=1,2,\ldots\}$. Let $\bar\W$ be the closure of $\W$
in $\h$. It readily follows from the properties
of $\breve\W$ that $\bar\W$ is compact and contains
(forward) paths $\eta_{x,t}$, not necessarily unique,
starting from each
$(x,t)\in\bar\R^2$, all of which are semipaths of $\breve\W$. Fix
now $i_0\geq0$ and let
\begin{eqnarray*}
{\mathbb D}^+\!\!\!&=&\{(x,t)\in\bar\R^2:\,t\leq t_{i_0},\,\eta_{x,t}(t_{i_0})\geq W_{i_0}(t_{i_0})
\mbox{ for some } \eta_{x,t}\in\bar\W\},\\
{\mathbb D}^-\!\!\!&=&\{(x,t)\in\bar\R^2:\,t\leq t_{i_0},\,\eta_{x,t}(t_{i_0})\leq W_{i_0}(t_{i_0})
\mbox{ for some } \eta_{x,t}\in\bar\W\},\\
{\mathbb D}&=&\,{\mathbb D}^+\cap{\mathbb D}^-.
\end{eqnarray*}
We claim that
${\mathbb D}=W^b_{i_0}$ (as a subset of $\bar\R^2$). Indeed, given
$(y,s)\in W^b_{i_0}$ with $s<t_{i_0}$, if $(y_n,s_n)_{n\geq1}$ approaches
$(y,s)$ with $y_n>y$ for $n\geq1$, then, by the compactness of $\bar\W$,
$(\eta_{y_n,s_n})_{n\geq1}$ has a limit
point $\eta_{y,s}$ in $\bar\W$ starting from $(y,s)$, which being a semipath
from a path of $\breve\W$ does not cross $W^b_{i_0}$. Thus
$(y,s)\in{\mathbb D}^+$. With an analogous argument
taking $y_n < y$, we conclude
that $(y,s)\in{\mathbb D}^-$, and thus $W^b_{i_0}\subset{\mathbb D}$.

If $(y,s)\in{\mathbb D}\setminus W^b_{i_0}$, then there is a semipath
$\eta_{y,s}$ of $\breve\W$ starting in $(y,s)$ which passes through
$(x_{i_0},t_{i_0})$, and thus a path in $\breve\W$ through
$(x_{i_0},t_{i_0})$ other than $\gamma_{x_{i_0},t_{i_0}}$, contradicting
the hypothesis. $\square$

\bigskip

\noindent{\bf Proof of Lemma~\ref{lem:l3}} It suffices to show that there is only one path in
$\breve\W$ through each point of ${\cal D}$, and then invoke
Lemma~\ref{lem:l2}. Given $(x_0,t_0)\in{\cal D}$, we have at least one path in $\breve\W$
through it, since we already have the semipath $W_{x_0,t_0}$. 
If there were two different
paths $\gamma$ and $\gamma'$ in $\breve\W$ through $(x_0,t_0)$, then, since
${\cal D}$ is dense, there would exist a point
$(x_1,t_1)$ in ${\cal D}$ with $t_1<t_0$ and $x_1$ strictly between $\gamma$ and
$\gamma'$, and thus, by the noncrossing property, $W_{x_1,t_1}$ would not be
able to avoid $(x_0,t_0)$, contradicting the hypothesis
that each $W_{x_i,t_i}$ avoids all other points of ${\cal D}$. $\square$

\bigskip

\noindent{\bf Proof of Proposition~\ref{prop:dbw}} Fix a
deterministic countable dense $\cd\in\R^2$ and let $\W=\W(\cd)$
be the set of paths of $\bar\W$ starting from $\cd$. Then $\W$
is almost surely a Brownian web skeleton. 
To fix notation, we write
$\cd=\{(x_i,t_i),\,i\geq1\}$, $\W=\{\hat W_i,\,i\geq1\}$, where
$\hat W_i$ is the path of $\W$ starting from $(x_i,t_i)$,
$i\geq1$. We proceed, as in~\cite{kn:STW, kn:FINR2}
and with the same notation, to 
define sets of coalescing/reflecting dual paths 
%follows (see Section 3 of~\cite{kn:FINR2}). Define $\hat
$\hat W_1^{b,n},\hat W_2^{b,n},\ldots$ inductively, as follows. Let $B_j^b$,
$j\geq1$ be i.i.d.~standard Brownian motions which are independent
of $\W$, and let $W_j^b$ be as in Equation~3.2 of~\cite{kn:FINR2}.
For $n>1$, let 
\begin{eqnarray}
 \label{eq:d1}
  \hat W_1^{b,n}&=&
  CR(W_1^b;\,\hat W_1,\hat W_2,\ldots,\hat W_{n});\\
   \label{eq:d2}
   \hat W_j^{b,n}&=&
   CR(W_j^b;\,\hat W_1,\ldots,\hat W_{n},\,\hat W_1^{b,n},\ldots,\hat
   W_{j-1}^{b,n}),\,1<j\leq n.
\end{eqnarray}
By Theorem 8 of~\cite{kn:STW},
$\{\hat W_1,\ldots,\hat W_{n},\,\hat W_1^{b,n},\ldots,\hat W_{j}^{b,n}\}$
has the same distribution as
$\{\tilde W_1,\ldots,\tilde W_{n},\,\tilde W_1^{b},\ldots,\tilde
W_{j}^{b}\}$ defined in Equations (3.3)-(3.5) of~\cite{kn:FINR2}. Since the
latter converges weakly to the DBW, then so does the former.

We claim that for each $i\geq1$, $\G_i:=[\{\hat
W_i^{b,n},\,1<n\leq m\},\,m>1]$ is tight. Thus by Proposition B.4
of~\cite{kn:FINR1}, $\{\hat W_i^{b,n},\,1<n<\infty\}$ is
a.s.~precompact.

We argue now that $\hat W_i^{b,n}\to\hat W_i^{b}$ almost surely as
$n\to\infty$ for all $i\geq1$. If that were not the case, by the
above precompactness, there would exist two different
subsequential limits for $\{\hat W_i^{b,n},\,1<n<\infty\}$, say
$\check W_i^{b}$ and $\breve W_i^{b}$, both backward paths
starting at $(x_i,t_i)$ and not crossing $\W$. There must thus
exist $(x_{i'},t_{i'})\in\cd$ with $t_{i'}<t_{i}$ and $x_{i'}$
strictly between $\check W_i^{b}(t_{i'})$ and $\breve
W_i^{b}(t_{i'})$. This and noncrossing imply that
$W_{i'}(t_i)=x_i$, which is a null event.

We thus have that for all $k\geq1$,
\begin{equation}\nn
  \{\hat W_1,\ldots,\hat W_{k},\,\hat W_1^{b,n},\ldots,\hat W_{k}^{b,n}\}\to
  \{\hat W_1,\ldots,\hat W_{k},\,\hat W_1^{b},\ldots,\hat W_{k}^{b}\}\mbox{ as
  } n\to\infty.
\end{equation}
By one of the above arguments, $\{\hat W_1,\ldots,\hat W_{k},\,\hat
W_1^{b},\ldots,\hat W_{k}^{b}\}$ has the same distribution as
$\{\tilde W_1,\ldots,\tilde W_{n},\,\tilde W_1^{b},\ldots,\tilde
W_{j}^{b}\}$ defined in Equations (3.3)-(3.5) of~\cite{kn:FINR2}.
It follows that
\begin{equation}
   \label{eq:d4}
  \{\hat W_1,\ldots,\hat W_{k},\,\hat W_1^{b},\ldots,\hat W_{k}^{b}\}
  \to(\hat\W,\hat\W^b)\mbox{ as } k\to\infty,
\end{equation}
where $(\hat\W,\hat\W^b)$ is a standard DBW
and the convergence is in $\h \times \h^b$. Clearly $\hat\W=\bar\W$.

It remains to justify the above tightness claim. We verify the
tightness condition $(T_1)$ from Appendix B of~\cite{kn:FINR1} for
$\G_i$. The reasoning goes by a {\em blocking argument} which is a
entirely similar to the one used for Proposition B.3 in that reference.
To facilitate the adaptation of that argument, let $\X^b_m=\{\hat
W_i^{b,n},\,1<n\leq m\}$ and $\X_m=\{\hat W_1,\ldots,\hat
W_{m}\}$. We then have that the paths of $\X^b_m$ do not cross
those of $\X_m$, and the latter satisfies condition $(I'_1)$ of
that proposition.

The blocking argument of the proof of Proposition B.3
in~\cite{kn:FINR1}
consists of subdividing the intervals $[x_0-u/2,x_0-u/4]\times\{t_0\}$ and
$[x_0+u/2,x_0+u/4]\times\{t_0\}$ into subintervals whose lengths are of order
$\sqrt t$, and then arguing that the $\limsup_{m\to\infty}$ of the probability
of either not having a forward crossing by any of the paths of $\X_m$ of any of the
rectangles to the right of $(x_0,t_0)$ with height $2t$ 
and having the respective
subintervals as base, or of not having a forward crossing by any of the paths of
$\X_m$ of any of the
rectangles to the left of $(x_0,t_0)$ with height $2t$ having the respective
subintervals as base is $o(t^{-1})$.
In the present case the same of course holds true, and as soon as there is a
blocking at each side, there cannot be a horizontal crossing by paths in
$\X^b_m$. Condition $(T_1)$ is thus verified and tightness of $\G_i$ follows.
$\square$

%%%%%%%%%%%%%%%%%%%%%%%%%%%%%%%%%%%%%%%%%%%%%%%%%%%%%%%%%%%%%%%%%%%%%%%%%%%%
%%%%%%%%%%%%%%%% REFERENCES %%%%%%%%%%%%%%%%%%%%%%%%%%%%%%%%%%%%%%%%%%%%%%%%%%%%%%
%%%%%%%%%%%%%%%%%%%%%%%%%%%%%%%%%%%%%%%%%%%%%%%%%%%%%%%%%%%%%%%%%%%%%%%%%%%%%

\end{document}